\documentclass{amsart}
\usepackage{amssymb,amsfonts,amsmath,amsthm,amscd,biblatex,stmaryrd,dsfont,esint,hyperref,upgreek,xcolor}
\usepackage[title]{appendix}
\usepackage[mathcal]{euscript}

\usepackage{fullpage}
\theoremstyle{plain}
\newtheorem{thm}{Theorem}
\newtheorem{cor}{Corollary}

\theoremstyle{definition}
\newtheorem{defn}[cor]{Definition}

\usepackage{biblatex}
{}
{}

\DeclareMathOperator{\sgn}{sgn}

\DeclareMathOperator{\dist}{dist}

\newcommand{\R}{\mathbb{R}}

\newcommand{\Z}{\mathbb{Z}}

\renewcommand{\tilde}{\widetilde}
\renewcommand{\hat}{\widehat}
\newcommand{\Minnie}{\textbf{\RN{1}}}
\newcommand{\Maxie}{\textbf{\RN{2}}}

\title{Slow periodic homogenization for Hamilton\textendash{}Jacobi equations}
\author{William Cooperman}
\begin{document}
\begin{abstract}
    Capuzzo-Dolcetta\textendash{}Ishii proved that the rate of periodic homogenization for coercive Hamilton\textendash{}Jacobi equations is $O(\varepsilon^{1/3})$. We complement this result by constructing examples of coercive nonconvex Hamiltonians whose rate of periodic homogenization is $\Omega(\varepsilon^{1/2})$.
\end{abstract}
\maketitle
\section{Introduction}
Since Lions\textendash{}Papanicolaou\textendash{}Varadhan~\cite{LPV} proved periodic homogenization of coercive Hamilton\textendash{}Jacobi equations, quantifying the rate of convergence has been a well-known open problem in both periodic and random settings. In a periodic environment, without additional structural assumptions on the Hamiltonian, the best known result so far is the $O(\varepsilon^{1/3})$ rate, proven by Capuzzo-Dolcetta\textendash{}Ishii~\cite{CDI}, which was also the first quantitive bound. On the other hand, when the Hamiltonian is convex in the momentum variable, the optimal rate of $O(\varepsilon)$ can be deduced from the optimal control formulation and an argument of Burago~\cite{Burago}, who proved a corresponding rate for homogenization of $\Z^d$-periodic metrics on $\R^d$. It is therefore natural to ask whether, in the absence of convexity, the $O(\varepsilon)$ rate still holds. Indeed, Ziliotto's~\cite{Ziliotto} example of stochastic non-homogenization (see also Feldman\textendash{}Souganidis~\cite{FeldSoug}) suggests that saddle points of the Hamiltonian may play a key role in slowing down periodic homogenization.

In this note, we answer this question in the negative by constructing examples which homogenize at a rate of $\Theta(\varepsilon^{1/2})$. In dimensions $d \geq 3$, we construct an example so that the effective Hamiltonian also happens to be convex.

Suppose that the Hamiltonian $H \colon \R^d \times \R^d \to \R$ is locally Lipschitz, $\Z^d$-periodic in the first variable, $x \in \R^d$, and uniformly coercive in the second variable, $p \in \R^d$; that is,
\[
    \liminf_{|p| \to \infty} \inf_{x \in \R^d} H(x, p) = +\infty.
\]

The microscopic problem at scale $\varepsilon > 0$ is the initial-value problem
\begin{equation}\label{micro-problem}
    \begin{cases}
        D_t u^\varepsilon(t, x) + H(\varepsilon^{-1}x, D_x u^\varepsilon(t, x)) = 0 &\quad \text{for $t > 0$ and $x \in \R^d$,}\\
        u^\varepsilon(0, x) = u_0(x) &\quad \text{for $x \in \R^d$,}
    \end{cases}
\end{equation}
where the initial data $u_0$ is Lipschitz.

Lions\textendash{}Papanicolaou\textendash{}Varadhan~\cite{LPV} proved that there is an effective Hamiltonian $\overline{H} \colon \R^d \to \R$, uniquely determined by $H$, such that $u^\varepsilon \to \overline{u}$ uniformly on compact sets, where $\overline{u}$ solves the effective problem
\begin{equation}\label{effective-problem}
    \begin{cases}
        D_t \overline{u}(t, x) + \overline{H}(D_x \overline{u}(t, x)) = 0 &\quad \text{for $t > 0$ and $x \in \R^d$,}\\
        \overline{u}(0, x) = u_0(x) &\quad \text{for $x \in \R^d$.}
    \end{cases}
\end{equation}

Under these general assumptions, the only known quantitative upper bound on the rate of homogenization is due to Capuzzo-Dolcetta\textendash{}Ishii~\cite{CDI}, who proved that $\|u^\varepsilon - \bar{u}\|_{L^\infty([0, T] \times \R^d)} = O(\varepsilon^{1/3})$ for $T > 0$. 

\begin{thm}\label{thm:main}
    There exists  a locally Lipschitz Hamiltonian $H \colon \R^d \times \R^d \to \R$, $\Z^d$-periodic in the first variable and uniformly coercive in the second variable, along with initial data $u_0 \in C^{0,1}(\R^d)$, such that, for sufficiently small $\varepsilon$,
    \[
        c\varepsilon^{1/2} \leq u^\varepsilon(1, 0) \leq C\varepsilon^{1/2}
    \]
    where $u^\varepsilon$ is the solution to the microscopic problem~\eqref{micro-problem}. If $d \geq 3$, then, furthermore, $H$ can be chosen so that $\overline{H}$ is convex.
\end{thm}

\section{Examples of slow homogenization}
To construct the examples, we first recall some facts from the theory of differential games. For a more thorough treatment, see Isaacs~\cite{Isaacs} and Evans\textendash{}Souganidis~\cite{EvansSoug}.

Let $A, B \subseteq \R^d$ be compact sets. We consider a differential game between two players named \Minnie{} and \Maxie{}. The game has a score, which \Minnie{} tries to minimize and \Maxie{} tries to maximize.
\begin{defn}
    A control for \Minnie{} (resp. \Maxie{}) is a measurable function $a \colon \R_{\geq 0} \to A$ (resp. $b \colon \R_{\geq 0} \to B$). We write $\mathcal{C}_A, \mathcal{C}_B$ to denote the set of controls for \Minnie{} and \Maxie{} respectively.
\end{defn}
\begin{defn}
    A strategy for \Minnie{} is a function $\alpha \colon \mathcal{C}_B \to \mathcal{C}_A$ with the nonanticipative property: if $t > 0$ and $b_1, b_2 \in \mathcal{C}_B$ with $b_1(s) = b_2(s)$ for almost all $s \in [0, t]$, then $\alpha(b_1)(s) = \alpha(b_2)(s)$ for almost all $s \in [0, t]$ also. We write $\mathcal{S}_A$ to denote the set of strategies for \Minnie{}, and define the set of strategies $\mathcal{S}_B$ for \Maxie{} correspondingly.
\end{defn}

A differential game is specified by the sets $A, B$, some Lipschitz initial data $u_0 \in C^{0,1}(\R^d)$, a running cost $R \in L^\infty(\R^d \times A \times B)$, and a transition function $f \in L^\infty(\R^d \times A \times B; \R^d)$ which is Carath\'eodory, i.e. $f(x, a, b)$ is continuous in $a, b$ for fixed $x$, and measurable in $x$ for fixed $a, b$. The game is based on the evolution of the state, $\sigma \colon \R_{\geq 0} \to \R^d$. Given a strategy $\alpha$, a control $b$, and a starting state $x$, the state evolves to satisfy the ordinary differential equation
\begin{equation}
    \begin{cases}
        \dot{\sigma}(t) = f(\sigma(t), \alpha(b)(t), b(t)) &\quad \text{for $t > 0$,}\\
        \sigma(0) = x &\quad \text{at $t=0$.}
    \end{cases}
\end{equation}
Sometimes, we will write $\sigma(t) = \sigma(t, \alpha, b)$ to emphasize the dependence on $\alpha$ and $b$.

Given $t \geq 0$ and $x \in \R^d$, the upper value of the game starting at $x$ after time $t$ is defined by
\begin{equation}\label{u-formula}
    u^+(t, x) := \inf_{\alpha \in \mathcal{S}_A} \sup_{b \in \mathcal{C}_B} \int_0^t R(\sigma(s), \alpha(b)(s), b(s)) \; ds + u_0(\sigma(t)).
\end{equation}

The upper value of the game (see Evans\textendash{}Souganidis~[\cite{EvansSoug}, Theorem 4.1] and Lions~\cite{Lions}) is the viscosity solution of the initial-value problem
\begin{equation}
    \begin{cases}
        D_t u^+(t, x) + H^+(x, D_x u^+(t, x)) = 0 &\quad \text{for $t > 0$ and $x \in \R^d$,}\\
        u^+(0, x) = u_0(x) &\quad \text{for $x \in \R^d$,}
    \end{cases}
\end{equation}
where the upper Hamiltonian $H \colon \R^d \times \R^d \to \R$ is given by
\begin{equation}\label{H-formula}
    H^+(x, p) := -\min_{a \in A} \max_{b \in B} R(x, a, b) + p \cdot f(x, a, b).
\end{equation}

It is worth noting that, by interchanging the order of the players, we can similarly define the lower value of the game by
\begin{equation}\label{u-minus-formula}
    u^-(t, x) := \sup_{\beta \in \mathcal{S}_B} \inf_{a \in \mathcal{C}_A} \int_0^t R(\sigma(s), a(s), \beta(a)(s)) \; ds + u_0(\sigma(t)),
\end{equation}
which solves a similar initial-value problem corresponding to the lower Hamiltonian
\begin{equation}\label{H-minus-formula}
    H^-(x, p) := -\max_{b \in B} \min_{a \in A} R(x, a, b) + p \cdot f(x, a, b).
\end{equation}

In all of our examples, the Isaacs condition $H^+ = H^-$ will be satisfied and therefore the upper and lower values of the game coincide. For brevity, we put $H := H^- = H^+$.

\subsection{An example in two dimensions}
Let $\varphi \colon \R/\Z \to [0, 1]$ be smooth such that $\varphi(0) = 1$, and $\varphi(x) = 0$ if $|x| \geq \frac{1}{100}$. We take $d=2$ and $A = \overline{B_1(0)}$, the closed unit ball centered at the origin, and $B = [0, 1] \times \{0\}$. Define the running cost by
\begin{equation}\label{ex1-running}
    R(x, a, b) := 100\left(1 - \varphi\left(x_2\right) - \varphi\left(x_2+\frac12\right)\right) + 100|a|^2,
\end{equation}
the transition function by
\begin{equation}\label{ex1-transition}
    f(x, a, b) := 2a + b\left(\varphi\left(x_2\right)-\varphi\left(x_2+\frac12\right)\right),
\end{equation}
and the initial data by
\[
    u_0(x) := \min\{|x_1|, 1\}.
\]

Although it's unnecessary for the proof, we note that the Isaacs condition
\[
    H(x, p) = -\min_{a \in A} \max_{b \in B} R(x, a, b) + p \cdot f(x, a, b) = -\max_{b \in B} \min_{a \in A} R(x, a, b) + p \cdot f(x, a, b)
\]
is satisfied, so the upper and lower values of this game coincide.

Intuitively, the microscopic environment consists of horizontal ``highways'' at every height in $\frac{\varepsilon}{2}\Z$. Outside of the highways, the running cost is punishingly large, so \Minnie{} is forced to spend most of the time inside the highways. Outside the highways, \Maxie{}'s control has no affect on the state. Inside highways at height in $\varepsilon\Z$, \Maxie{} has the option to push the state in the $+e_1$ direction, and inside highways at height in $\varepsilon\left(\Z+\frac12\right)$, \Maxie{} has the option to push the state in the $-e_1$ direction. \Maxie{}'s control has no effect on the running cost, and \Minnie{} is heavily penalized for pushing the state in any direction. If \Minnie{} wants to stay close to the origin (where the terminal cost is lowest), then one strategy is to enter a highway and wait until \Maxie{} pushes the state a distance of $\varepsilon^{1/2}$ from the origin. Then, \Minnie{} can switch to a highway that leads back to the origin, and repeat. By the same reasoning, \Maxie{} can force \Minnie{} to switch highways at least $\varepsilon^{-1/2}$ many times, or else \Minnie{} risks paying a terminal cost of at least $\varepsilon^{1/2}$. Each highway switch adds running cost proportional to $\varepsilon$ to the total, so the error terms balance.

Now, we prove the $d=2$ case of Theorem~\ref{thm:main}.
\begin{proof}
    We use the differential game characterization~\eqref{u-formula} of $u^\varepsilon$. For the upper bound, we construct a strategy $\alpha \colon \mathcal{C}_B \to \mathcal{C}_A$ for \Minnie{} piecewise as follows.
    \begin{enumerate}
        \item At the beginning of this step, suppose that the strategy has already been constructed up to time $t_i \geq 0$ and $\sigma(t_i, \alpha, b) \in [-\frac{\varepsilon}{4}, \frac{\varepsilon}{4}] \times \{0\}$. Given $b \in \mathcal{C}_B$, let \[ t_{i+1} := \min \{ t > 0 \mid u_0(\sigma(t, \hat{\alpha}_{t_i}, b)) \geq \varepsilon^{1/2}\}, \] where
            \[
                \hat{\alpha}_{s}(t) := 
                \begin{cases}
                    \alpha(t) &\quad \text{for $t < s$,}\\
                    0 &\quad \text{for $t \geq s$.}
                \end{cases}
            \]
            From the structure of $f$ and $u_0$, we deduce that $\sigma(t_{i+1}, \hat{\alpha}_{t_i}, b) = (\varepsilon^{1/2}, 0)$. Set $\alpha(b)(t) := \hat{\alpha}_{t_i}(b)(t)$ for $t < t_{i+1}$.
        \item Write $t_{i+2} := t_{i+1}+\frac{\varepsilon}{4}$ and set $\alpha(b)(t) := (0, 1)$ for $t_{i+1} \leq t < t_{i+2}$.
        \item We deduce that $\sigma(t_{i+2}, \alpha, b) \in \left[\varepsilon^{1/2} - \frac{\varepsilon}{4}, \varepsilon^{1/2} + \frac{\varepsilon}{4}\right] \times \left\{\frac{\varepsilon}{2}\right\}$. Now, set \[ t_{i+3} := \min \{ t > 0 \mid u_0(\sigma(t, \hat{\alpha}_{t_{i+2}}, b)) \leq 0\} \] and set $\alpha(b)(t) := \hat{\alpha}_{t_{i+2}}(b)(t)$ for $t_{i+2}  \leq t < t_{i+3}$. As in the first step, we deduce that $\sigma(t_{i+3}, \alpha, b)  = \left(0, \frac{\varepsilon}{2}\right)$.
        \item Write $t_{i+4} := t_{i+3}+\frac{\varepsilon}{4}$ and set $\alpha(b)(t) := (0, -1)$ for $t_{i+3} \leq t < t_{i+4}$.
        \item We deduce that $\sigma(t_{i+4}, \alpha, b) \in \left[-\frac{\varepsilon}{4}, \frac{\varepsilon}{4}\right] \times \{0\}$. Now, go back to step 1 and repeat, but starting at time $t_{i+4}$ instead of time $t_i$.
    \end{enumerate}
    The construction maintains the invariant that $\sigma(t, \alpha, b) \in \left[0, \frac{\varepsilon}{4}\right] \times \left[-\frac{\varepsilon}{4}, \varepsilon^{1/2} + \frac{\varepsilon}{4}\right]$, so we ensure that the terminal cost is at most $u_0(\sigma(1)) \leq \varepsilon^{1/2} + \frac{\varepsilon}{4} \leq 2\varepsilon^{1/2}$.

    It remains to show that the running cost given by $\alpha$ is at most $C\varepsilon^{1/2}$. For any $i=0,1,2,\dots$, we claim
    \[
        \int_{t_i}^{t_{i+1}} R(\varepsilon^{-1}\sigma(t), \alpha(b)(t), b(t)) \; \mathrm{d}t \leq 50\varepsilon.
    \]
    Indeed, an interval created in step 2 or step 4 above satisfies this bound, as $t_{i+1}-t_i = \frac{\varepsilon}{4}$ and $R \leq 200$ everywhere. On the other hand, intervals created in step 1 or step 3 above have running cost $0$, since $\sigma(t) \in \R \times \frac{\varepsilon}{2}\Z$ for all $t$ in the interval.

    We have shown that each step adds at most $50\varepsilon$ to the running cost. On the other hand, every interval created by step 1 or step 3 runs for time at least $\varepsilon^{1/2}$, so there can be at most $\varepsilon^{-1/2}$ such intervals in $[0, 1]$. We conclude that the total running cost is at most $50\varepsilon^{1/2}$, so $C=50$ satisfies the claim.

    Next, we turn to the lower bound $u^\varepsilon(1, 0) \geq c\varepsilon^{1/2}$. Given a strategy $\alpha \colon \mathcal{C}_B \to \mathcal{C}_A$, we construct the following control for \Maxie{}.
    \[
        b(t) := \begin{cases} (1, 0) &\quad \text{if $u_0(\sigma(t, \alpha, b)) \leq 2\varepsilon^{1/2}$ or $\sigma(t, \alpha, b) \cdot f(\varepsilon^{-1}\sigma(t, \alpha, b), 0, b) > 0$,}\\ 0 &\quad \text{otherwise.} \end{cases}
    \]
    We claim that this control yields a value of at least $c\varepsilon^{1/2}$. Indeed, consider the set of times
    \[
        E := \{t \in [0, 1] \mid R(\varepsilon^{-1}\sigma(t, \alpha, b), 0, 0) \geq 1\}.
    \]
    First, note that
    \begin{equation}\label{E-bound}
        \int_0^1 R(\varepsilon^{-1}\sigma(t, \alpha, b), \alpha(b)(t), b(t)) \; \mathrm{d}t \geq |E|,
    \end{equation}
    where $|E|$ denotes the Lebesgue measure. On the other hand, let $U_+ := \{t \in [0, 1] \setminus E \mid \varphi(\varepsilon^{-1}{\sigma(t, \alpha, b)}_2) > 0\}$ and $U_- := [0, 1] \setminus (E \cup U_+)$. In the set of times $U_+$ (resp. $U_-$), \Maxie{} can control the state to push in the $+e_1$ (resp. $-e_1$) direction with magnitude at least $0.99$.

    There are two cases.
    \begin{enumerate}
        \item Suppose that there exist $0 \leq t_0 \leq t_1 \leq 1$ such that $|[t_0, t_1] \cap U_-| \geq 6\varepsilon^{1/2}$ and $|[t_0, t_1] \cap U_+| = 0$. Then either
            \[
                \int_{t_0}^{t_1} R(\sigma(t, \alpha, b), \alpha(b)(t), b(t)) \; \mathrm{d}t \geq \varepsilon^{1/2},
            \] in which case we conclude, or we have $u_0(\sigma(t_1, \alpha, b)) \geq 3\varepsilon^{1/2}$, because \Minnie{} has a total effect of less than $\frac{\varepsilon^{1/2}}{100}$ on the state over the interval $[t_0, t_1]$. If $t$ is such that $u_0(\sigma(t, \alpha, b)) \geq 2\varepsilon^{1/2}$, then by definition of $b$ we have $\sigma(t, \alpha, b) \cdot f(\varepsilon^{-1}\sigma(t, \alpha, b), 0, b) \geq 0$, so the control of \Maxie{} never pushes the state toward smaller values of $u_0$. So, either $u(\sigma(1, \alpha, b)) \geq 2\varepsilon^{1/2}$ or \Maxie{} spends
            \[
                \int_{t_l}^{1} R(\sigma(t, \alpha, b), \alpha(b)(t), b(t)) \; \mathrm{d}t \geq \varepsilon^{1/2},
            \]
            and in either case we conclude. We note that if $|[t_0, t_1] \cap U_+| \geq 6\varepsilon^{1/2}$ and $|[t_0, t_1] \cap U_-| = 0$, then we conclude by the same argument.
        \item Otherwise, we may assume by~\eqref{E-bound} that $|E| \leq \varepsilon^{1/2}$, so $|U_1 \cup U_2| \geq 1-\varepsilon^{1/2}$. Take $t_- \in U_-$ and $t_+ \in U_+$. Since $|{\sigma(t_-, \alpha, b)}_2-{\sigma(t_+, \alpha, b)}_2| \geq \varepsilon\left(\frac{1}{2}-\frac{1}{50}\right)$ (using the fact that $\varphi$ is supported in $\left(-\frac{1}{100}, \frac{1}{100}\right)$), we conclude that $|[t_-, t_+] \cap E| \geq \frac{\varepsilon}{5}$ and therefore
            \[
                \int_{t_-}^{t_+} R(\sigma(t, \alpha, b), \alpha(b)(t), b(t)) \; \mathrm{d}t \geq \frac{\varepsilon}{5}.
            \]
            On the other hand, since the hypotheses of the previous case do not apply, we have $[t_0, t_1] \cap U_- \neq \emptyset$ and $[t_0, t_1] \cap U_+ \neq \emptyset$ whenever $t_1-t_0 \geq 6\varepsilon^{1/2}$. So, there are at least $\frac{\varepsilon^{-1/2}}{6}-1$ many such disjoint intervals in $[0, 1]$, and we conclude that
            \[
                |E| \geq \frac{\varepsilon}{5} \cdot \left(\frac{\varepsilon^{-1/2}}{6}-1\right) \geq \frac{\varepsilon^{1/2}}{35}
            \]
            as long as $\varepsilon \leq \frac{1}{42^2}$, and therefore
            \[
                \int_{0}^{1} R(\sigma(t, \alpha, b), \alpha(b)(t), b(t)) \; \mathrm{d}t \geq \frac{\varepsilon^{1/2}}{35}.
            \]
    \end{enumerate}
    In any of the cases, $c=\frac{1}{35}$ satisfies the claim.
\end{proof}
\subsection{An example in three and higher dimensions}
Next, we show that if $d \geq 3$, we can construct an example where the effective Hamiltonian $\overline{H}$ is convex. Without loss of generality, let $d=3$. Inspired by the example of Hedlund~\cite{Hedlund}, let
\[
    \mathcal{L} := \bigcup_{i=1}^3 \ell_i + \Z^3,
\]
where $\ell_1 := \R \times \{0\} \times \{0\}$, $\ell_2 := \{0\} \times \R \times \{\frac14\}$, and $\ell_3 := \{\frac14\} \times \{\frac14\} \times \R$, and for $i \in \{1,2,3\}$ let $\varphi_i \colon \R^3/\Z^3 \to [0, 1]$ be smooth such that $\varphi_i(x) = 1$ for $x \in \ell_i$ and $\varphi_i(x) = 0$ if $\dist(x, \ell_i) \geq \frac{1}{100}$. Write $\varphi := \varphi_1 + \varphi_2 + \varphi_3$ and $\tilde{\varphi} := (\varphi_1, \varphi_2, \varphi_3)$.

Let $A = \overline{B_2(0)}$ and $B = {[0, 1]}^3$. Define the running cost by
\begin{equation}\label{ex2-running}
    R(x, a, b) := 100\left(1 - \varphi\left(x\right) - \varphi\left(x+\frac12\right)\right) + 100|a|,
\end{equation}
and the transition function by
\begin{equation}\label{ex2-transition}
    f(x, a, b) := 2\left(1+99\left(\varphi\left(x\right)+\varphi\left(x+\frac12\right)\right)\right)a + b \odot \left(\tilde{\varphi}\left(x\right)-\tilde{\varphi}\left(x+\frac12\right)\right),
\end{equation}
where we write $\frac12$ to denote the vector $(\frac12, \frac12, \frac12)$ and $\odot$ to denote the pointwise product, i.e. $x \odot y = (x_1y_1, x_2y_2, x_3y_3)$.

As before, the Isaacs condition is satisfied and therefore the upper and lower values of the game coincide. This environment is very similar to the previous example, except that now the highways go in every coordinate direction $\pm e_i$. We take advantage of the fact that this is possible in three dimensions while ensuring that each highway is far away from any other highway. The running cost is nearly identical to the previous example, except that the penalty for \Minnie{} is $100|a|$ instead of $100|a|^2$. The transition function is slightly different; \Maxie{}'s controls work similarly, but \Minnie{} is now able to move faster inside the highways and slower outside.

Now, we prove the $d \geq 3$ case of Theorem~\ref{thm:main}.
\begin{proof}
    Using the same initial data $u_0(x) := \min\{|x_1|, 1\}$, the argument for the $\varepsilon^{1/2}$ rate is identical to the previous example, so we omit it. To show that $\overline{H}$ is convex, we find, for each $p \in \R^d$ and $\lambda > 0$, bounds for the long-time corrector $v^{p}$, defined as the solution to the initial-value problem
    \begin{equation}\label{approx-corrector-problem}
        \begin{cases}
            D_t v^{p}(t, x) + H(x, D_x v^{p}(t, x)) = 0 &\quad \text{for $t > 0$ and $x \in \R^2$,}\\
            v^{p}(0, x) = p \cdot x &\quad \text{for $x \in \R^2$.}
        \end{cases}
    \end{equation}
    Then, using the fact that $\overline{H}(p) = -\lim_{t \to \infty} t^{-1}v^p(t, 0)$, we obtain a formula for $\overline{H}$.

    First, for $\gamma \in \R$ we claim that
    \[
        \overline{H}(\gamma e_i) = h(\gamma) := \max \{0, 400|\gamma|-200\}.
    \]
    Indeed, we immediately have the lower bound
    \[
        v^{\gamma e_i}(t, x) \geq \min_{a \in A} t(100|a| + 200\gamma a \cdot e_i) = th(\gamma),
    \]
    which follows from considering the constant $0$ control for \Maxie{} and ignoring the space-dependent part of the running cost.

    On the other hand, to obtain the upper bound, \Minnie{} can use the following strategy: immediately move the state (in constant time) into $\ell_i+\Z^3$, if $\gamma < 0$, or into $\ell_i + \frac12 + \Z^3$ otherwise. For the rest of time, use the constant strategy $\alpha = -(\sgn \gamma)\min\left(2, |\gamma|\right)e_i$. Computing the result of the game with this strategy shows that
    \[
        v^{\gamma e_i}(t, x) \leq C + \min_{a \in A} t(100|a| + 200\gamma a \cdot e_i) = C + th(\gamma).
    \]

    We have shown that $\overline{H}(\gamma e_i) = h(\gamma)$, and $h(\gamma)$ is convex. To conclude, we will show that, for each $p \in \R^d$,
    \[
        \overline{H}(p) = \max_{i=1}^3 \overline{H}(p_i e_i).
    \]
    The inequality $\overline{H}(p) \geq \overline{H}(p_i e_i)$ follows immediately from using the strategy outlined above for $i \in \{1, 2, 3\}$.

    For the other inequality, let $\alpha \colon \mathcal{C}_B \to \mathcal{C}_A$ be a strategy for the initial data $p \cdot x$ and starting state $0$. Suppose for contradiction that for all $t > 0$ large, we have
    \[
        \sup_{b \in \mathcal{C}_B} \int_0^t R(\sigma(s, \alpha, b), \alpha(b)(s), b(s)) \; \mathrm{d}s + p \cdot \sigma(t, \alpha, b) \leq -t\max_{i=1}^3 h(p_i) - ct,
    \]
    for some small $c > 0$. We claim that this cannot hold even for the constant control $b = 0$. First, we note that the terminal cost can be interpreted as a kind of running cost, in the sense that
    \[
        p \cdot \sigma(t, \alpha, b) = \int_0^t p \cdot f(\sigma(s, \alpha, b), \alpha(b)(s), b(s)) \; \mathrm{d}s,
    \]
    since the starting state is $\sigma(0) = 0$.
    Let $E := \{t > 0 \mid R(\sigma(t, \alpha, 0), 0, 0) \geq 99\}$. Then if $t \in E$,
    \begin{align*}
        J(s) &:= R(\sigma(s, \alpha, b), \alpha(b)(s), b(s)) + p \cdot f(\sigma(s, \alpha, b), \alpha(b)(s), b(s))\\
        &\geq 99 + \min_{a \in A} \left[ 100|a| + 4 a \cdot p \right]\\
        &\geq 99 - \frac{1}{50}h(|p|)\\
        &\geq 99 - \frac{\sqrt{3}}{50}\max_{i=1}^3 h(p_i) - 4(\sqrt{3}-1)\\
        &\geq 95 - \frac{\sqrt{3}}{50}\max_{i=1}^3 h(p_i),
    \end{align*}
    where we use the fact that \Minnie{} can only move at much slower the speed away from the highways.

    On the other hand, we say that an interval $[t_0, t_1]$ stays close to a highway if there is a line $\ell$ in $\mathcal{L}$ such that, for every $s \in [t_0, t_1]$, the line $\ell$ is the closest line in $\mathcal{L}$ to $\sigma(s, \alpha, b)$. In any such interval $[t_0, t_1]$, we have
    \[
        \int_{t_0}^{t_1} J(s) \; \mathrm{d}s \geq -|p| - (t_1-t_0)\max_{i=1}^3 h(p_i),
    \]
    where the first term accounts for the (constant-sized) movement in the direction orthogonal to $\ell$, and the second term accounts for the movement in the direction parallel to $\ell$.

    Write $U := [0, t] \setminus E$. Given $t_0, t_1 \in U$, we write $t_0 \sim t_1$ iff $[t_0, t_1]$ stays close to a highway. Let $I_1, I_2, \dots, I_n$ denote the equivalence classes of $U/\sim$, ordered by the usual order on $\R$. We claim that $n \leq 40t$. Indeed, $[\max I_i, \min I_{i+1}] \subseteq E$, and $\min I_{i+1}-\max I_i \geq \frac{1}{40}$, since the distance $\frac{1}{100}$ neighborhoods of highways are at least distance $\frac{1}{5}$ apart, and the speed limit in $E$ is at most $8$.

    Putting everything together, we write
    \begin{align*}
        &\int_{0}^{t} R(\sigma(s, \alpha, b), \alpha(b)(s), b(s)) + p \cdot f(\sigma(s, \alpha, b), \alpha(b)(s), b(s)) \; \mathrm{d}s\\
        &\qquad \qquad = \int_0^{\min I_1} J(s) \; \mathrm{d}s + \sum_{j=1}^{n-1} \left(\int_{\min I_j}^{\max I_j} J(s) \; \mathrm{d}s + \int_{\max I_j}^{\min I_{j+1}} J(s) \; \mathrm{d}s\right)\\
        &\qquad \qquad \qquad + \int_{\min I_n}^{\max I_n} J(s) \; \mathrm{d}s + \int_{\max I_n}^t J(s) \; \mathrm{d}s\\
        &\qquad \qquad \geq -\frac{\sqrt{3}\min I_1}{50}\max_{i=1}^3 h(p_i) + 95(\min I_1)\\
        &\qquad \qquad \qquad + \sum_{j=1}^{n-1} - |p| - |I_j|\max_{i=1}^3 h(p_i) - \frac{\sqrt{3}(\min I_{j+1}-\max I_j)}{50}\max_{i=1}^3 h(p_i) + 95(\min I_{j+1}-\max I_j)\\
        &\qquad \qquad \qquad - |p| - |I_n|\max_{i=1}^3 h(p_i) - \frac{\sqrt{3}(t-\max I_n)}{50}\max_{i=1}^3 h(p_i)\\
        &\qquad \qquad \geq - t \max_{i=1}^3 h(p_i),
    \end{align*}
    where the sum telescopes and we use the fact that
    \[
        2|p| \leq \frac{1}{40}\left(95 + \left(1-\frac{\sqrt{3}}{50}\right)\max_{i=1}^3 h(p_i)\right).
    \]
\end{proof}
\section*{Acknowledgement} I would like to thank Charles Smart and Ahmed Bou-Rabee for many helpful discussions and comments.
\printbibliography{}
\end{document}